# Eigenvector-based identification of bipartite subgraphs

Debdas Paul[a,∗], Dragan Stevanović[b]

[a]*Institute for Systems Theory and Automatic Control, Pfaffenwaldring 9, 70569 Stuttgart, Germany*
[b]*Mathematical Institute of the Serbian Academy of Sciences and Arts, Knez Mihailova 36, 11001 Belgrade, Serbia*

**Abstract**

We report our experiments in identifying large bipartite subgraphs of simple connected graphs which are based on the sign pattern of eigenvectors belonging to the extremal eigenvalues of different graph matrices: adjacency, signless Laplacian, Laplacian, and normalized Laplacian matrix. We compare the performance of these methods to a 'local switching' algorithm based on the Erdös' bound that each graph contains a bipartite subgraph with at least half of its edges. Experiments with one scale-free and three random graph models, which cover a wide range of real-world networks, show that the methods based on the eigenvectors of the normalized Laplacian and the adjacency matrix yield slightly better, but comparable results to the local switching algorithm. We also formulate two edge bipartivity indices based on the former eigenvectors, and observe that the method of iterative removal of edges with maximum bipartivity index until one obtains a bipartite subgraph, yields comparable results to the local switching algorithm, and significantly better results than an analogous method that employs the edge bipartivity index of Estrada and Gomez-Gardeñes.

## 1. Introduction

A graph $G(V, E)$ with a set of vertices $V$ and a set of edges $E$ is bipartite if there exists a partition $V = X \cup Y$, $X \cap Y = \emptyset$ such that every edge $e \in E$ has one end in $X$ and another in $Y$. It is a classical result that a graph is bipartite if and only if it does not contain a cycle of odd length as a subgraph (Asratian et al., 1998).

Bipartite graphs have many important applications in various fields of science and technology. For example, in modern coding theory, bipartite graphs such as Factor graphs and Tanner graphs are used to decode codewords received from the channel (Moon, 2005). It has also been found that the underlying structure of many complex networks (Watts and Strogatz, 1998; Barabási and Albert, 1999; Guillaume and Latapy, 2004, 2006), like biological networks (Pavlopoulos et al., 2011, 2018; Platig et al., 2016; Baker et al., 2014; Moon et al., 2016), social networks (Newman et al., 2002; Schweitzer et al., 2009), and technological networks is bipartite in nature. In the emerging field of human-robot collaborative system design (Roy and Edan, 2018), bipartite graphs are used to improve the temporal coordination in human-robot teamwork (Chao and Thomaz, 2012). For a more comprehensive review of various applications of bipartite graphs, interested readers are referred to Holme et al. (2003). Recently bipartivity, a quantitative descriptor of bipartiteness of a graph, has been shown to be inversely correlated with the efficiency of airline transportation in Europe (Estrada and Gómez-Gardeñes, 2016).

Here we propose simple eigenvector-based methods aimed at identifying large bipartite subgraphs. The problem of identifying the largest bipartite subgraph of a graph is an instance of the MAX-CUT problem in which each edge has the same weight (Michael and David, 1979). From the computability point of view, the MAX-CUT problem is NP-hard (Michael and David, 1979; Goemans and Williamson, 1994) and there

---

∗Corresponding author
*Email addresses:* `debdas.paul@ist.uni-stuttgart.de` (Debdas Paul), `dragance106@yahoo.com` (Dragan Stevanović)

exists several approximation algorithms whose approximation ratio ranges from 0.5 (Mitzenmacher and Upfal, 2005; Motwani, 1995) to 0.879 (Goemans and Williamson, 1994) up to 0.942 (if P $\neq$ NP) (Håstad, 2001). The simplest approximation algorithm with an approximation ratio 0.5 is based on a theorem due to Erdös (Erdös, 1965), which shows that any graph $G(V, E)$ contains a bipartite subgraph with at least $\frac{|E|}{2}$ edges. This lower bound was later improved by Edwards to $\frac{|E|}{2} + \frac{|V|-1}{4}$ (Edwards, 1973). We refer to the method used in the proof of the Erdös' result as a *local switching* algorithm and implemented it in a programmatic manner for comparison purposes. Previously, Bylka et al. (1999) have introduced a local switching algorithm that guarantees the bound due to Edwards.

Motivation to consider eigenvector-based methods comes from a classical result in spectral graph theory which claims that a connected graph is bipartite if and only if for the spectral radius $\lambda_1$ of the adjacency matrix $A$ of $G$, $-\lambda_1$ is also an eigenvalue of $A$, in which case the eigenvector of $-\lambda_1$ is obtained from the principal eigenvector by changing signs of the components in one part of the bipartition (Sachs, 1964). A similar result holds for the normalized Laplacian matrix $\mathcal{L} = D^{-1/2}LD^{1/2}$, where $D$ is the diagonal matrix of vertex degrees, which claims that a connected graph is bipartite if and only if the largest eigenvalue of $\mathcal{L}$ is equal to 2 in which case the components of the eigenvector of 2 have equal value in one part and the opposite value in the other part of the bipartition (Chung, 1996). It is thus natural to expect that, in cases when a graph is close to being bipartite, the sign patterns of the eigenvector corresponding to the smallest eigenvalue of $A$ and the eigenvector corresponding to the largest eigenvalue of $\mathcal{L}$ will produce bipartitions of the vertex set containing most edges of the graph. The relation between bipartiteness and the spectral properties of $A$ had been further extended: if $B$ is an essentially non-negative symmetric matrix representing a connected, bipartite graph $G$ with bipartition $X \cup Y$, then there exists a non-zero eigenvector $z$ that belongs to the minimum eigenvalue of $B$ such that $z_i > 0$ for $i \in X$ and $z_j < 0$ for $j \in Y$ (Roth, 1989). In addition to $A$, the signless Laplacian matrix $Q = D + A$ is also *essentially* non-negative symmetric matrix, which partially justifies the use of the sign pattern of the eigenvector of the smallest eigenvalue of $Q$ as well. Further, it is known that the spectra of $Q$ and the Laplacian matrix $L = D - A$ coincide for a bipartite graph (Grone et al., 1990, prop. 2.2). If $G$ is a connected bipartite graph with bipartition $X \cup Y$, let $R$ be the diagonal matrix indexed by vertices of $G$ such that $r_{j,j} = 1$ if $j \in X$ and $r_{j,j} = -1$ if $j \in Y$. Observe that $R = R^{-1}$. It is straightforward to see that $R^{-1}LR = Q$, which yields a relationship between the eigenvectors of $Q$ and $L$: $Qz = \lambda z \iff L(Rz) = \lambda(Rz)$. The eigenvector of the smallest eigenvalue 0 of $L$ is the all-one vector $\mathbf{1}$ so the eigenvector of the smallest eigenvalue of $Q$ is $R\mathbf{1}$; note that $R\mathbf{1}$ has entries of one sign on $S$ and the opposite sign on $T$. The eigenvector of the largest eigenvalue of $Q$ is a positive Perron vector $v$, so that the eigenvector of the largest eigenvalue of $L$ is $Rv$; note that $Rv$ has entries of one sign on $S$ and the opposite sign on $T$, which partially justifies the use of the sign pattern of the eigenvector of the largest eigenvalue of $L$ as well.

To start with, we compare methods based on these particular eigenvectors to the local switching algorithm for three random graph models and one scale-free graph model, which cover a wide range of real-world networks. We observe that in these cases, methods based on the eigenvectors of $A$ and $\mathcal{L}$ yield slightly better, but comparable results to the local switching algorithm.

Further, Estrada and Rodríguez-Velázquez provided in (Estrada and Rodríguez-Velázquez, 2005) a measure of edge bipartivity that quantifies the contribution of edges toward bipartiteness based on the ratio of the number of closed even walks and the number of all closed walks present in a graph. The measure was further improved by Estrada and Gomez-Gardeñes (Estrada and Gómez-Gardeñes, 2016). We define analogous measures for edge bipartivity based on the eigenvectors that belong to the smallest eigenvalue of $A$ and the largest eigenvalues of $\mathcal{L}$, respectively, and then compare the performance of using these two measures in identifying large bipartite subgraphs. Simulation results show that these eigenvector-based measures outperform the analogous method that uses the Estrada and Gomez-Gardeñes' edge bipartivity.

## 2. Methods

We describe in this section algorithmic methods for identification of bipartite subgraphs. Notation used here and in subsequent sections is given in Table 1.

**Table 1. Notation**

| Notation | Meaning |
|---|---|
| $G(V, E)$ or $G$ | A simple, connected, unweighted, nonbipartite graph (unless explicitly stated otherwise) with vertex set $V$ and edge set $E$ |
| $A$, $L$, $Q$, $\mathcal{L}$, and $D$ | Adjacency, Laplacian, signless Laplacian, and normalized Laplacian, and degree matrix of $G$, respectively |
| $|E_{\text{bipart}}(G)|$ | Number of edges remaining when edge-deleted subgraph of $G$ becomes bipartite |
| $|E(G)|$ | Number of edges in $G$ |
| $r^b$ | $\frac{|E_{\text{bipart}}(G)|}{|E(G)|}$ |
| $|E_u^{\text{ext}}(G)|$ | Number of neighbors of vertex $u$ in the other part of the bipartition |
| $|E_u^{\text{int}}(G)|$ | Number of neighbors of vertex $u$ in its part of the bipartition |
| $\{\lambda_i^M\}$, $\lambda_i^M$ | The set of eigenvalues, and the $i^{\text{th}}$ largest eigenvalue of matrix $M$ |
| $\nu^{\lambda_i^M}$ | Eigenvector of the eigenvalue $\lambda_i^M$ of matrix $M$ |
| $\nu_k^{\lambda_i^M}$ | The $k^{\text{th}}$ entry of the eigenvector $\nu^{\lambda_i^M}$ of $M$ |
| $X \xrightarrow{i} Y$ | Movement of vertex $i$ from the set $X$ to $Y$ |
| E-R, W-S, RG, and B-A | Erdös-Rényi (Erdös and Rényi, 1959, 1960; Newman, 2003), Watts-Strogatz (Watts and Strogatz, 1998), Random Geometric (Gilbert, 1961; Penrose, 2003; Dall and Christensen, 2002; Pržulj et al., 2004), and Barabasi-Albert graph model(Barabási and Albert, 1999) respectively |

## 2.1. Local switching algorithm

The local switching algorithm proceeds by first randomly permuting the vertices and then partitioning them into sets $X$ and $Y$ ($X \cup Y = V$, $X \cap Y = \emptyset$) with sizes $\left\lfloor \frac{|V|}{2} \right\rfloor$ and $\left\lceil \frac{|V|}{2} \right\rceil$, respectively. Afterwards, a vertex $u$ is moved to another part if $|E_u^{\text{int}}(G)| > |E_u^{\text{ext}}(G)|$, or if $2 \times |E_u^{\text{int}}(G)| > |E_u(G)|$ since $|E_u^{\text{int}}(G)| + |E_u^{\text{ext}}(G)| = |E_u(G)|$. Vertex $u$ is then marked as *moved* and not considered for movement in subsequent iterations. The process continues alternately between the two parts until no more movements of vertices are possible. After the termination of the process, the ratio $r^b$ is calculated. Algorithm 1 summarizes the pseudo-code for the local switching algorithm.

**Remark 2.0.1.** *The ratio $r^b$ depends on the initial partition. The effect of the initial partition on $r^b$ is demonstrated via Example 2.1. In order to alleviate the problem, the maximum value of $r^b$ is taken over a number of different partitions.*

---

**Algorithm 1:** local switching algorithm

**Input**: $G = (V, E)$
**Result**: $r^b$
**Initialize**: $V_{\text{moved}} = \emptyset$; Partition randomly $X \cup Y = V \wedge X \cap Y = \emptyset$; Choose randomly $P = X, Q = Y$ or $P = Y, Q = X$; $T_{\text{Movement}} = 0$; $T_{\text{noMovement}} = 1$

1 **repeat**
2     **for** all $u \in P$ **do**
3         Compute $E_u^X$ and $E_u^Y$;
4         **if** $2E_u^P > E_u^V \wedge u \notin V_{\text{moved}}$ **then**
5             $Q \leftarrow u$;
6             $V_{\text{moved}} \leftarrow u$;
7             $T_{\text{Movement}} = 1$;
8             break;
9         **end**
10    **end**
11    **if** $T_{\text{Movement}} = 0$ **then**
12         $T_{\text{noMovement}} = T_{\text{noMovement}} + 1$;
13    **else**
14         $T_{\text{Movement}} = 0$;
15         $T_{\text{noMovement}} = 1$;
16    **end**
17    **if** $P = X$ **then**
18         $P = Y; Q = X$;
19    **else**
20         $P = X; Q = Y$;
21    **end**
22 **until** $T_{\text{noMovement}} > 2$;

**Example 2.1.** *Consider a simple, undirected, and connected graph $G(V, E)$ with $V = \{1, 2, 3, 4, 5\}$ and $E = \{(1, 2), (1, 3), (1, 5), (3, 4), (3, 5), (4, 5)\}$. Consider the following initial partition of $V$: $X = \{1, 2\}$ and $Y = \{3, 4, 5\}$. Initially*

$$E_1^X = 1, E_2^X = 1, E_3^Y = 3, E_4^Y = 2, E_5^Y = 2,$$

*while*

$$E_1^V = 3, E_2^V = 1, E_3^V = 3, E_4^V = 2, E_5^V = 2.$$

*If we start with $P = X$, the process continues as follows:*

1. *As $2E_i^X > E_i^V$ for $i = 2$, therefore $X \xrightarrow{2} Y$, $X = X \setminus \{2\}$ and $Y = Y \cup \{2\}$.*
2. *Now $2E_i^Y > E_i^V$ is satisfied for vertices 3, 4 and 5, so that any vertex among them is suitable for movement.*
3. *Let us pick vertex 3: $Y \xrightarrow{3} X$, $Y = Y \setminus \{3\}$ and $X = X \cup \{3\}$.*
4. *In the part $X = \{1, 3\}$, no movement of vertices from $X$ to $Y$ is possible as 3 is already been moved and 1 does not satisfy $2E_1^X > E_1^V$. Therefore $T_{\text{noMovement}} = 1$.*
5. *In the part $Y = \{2, 4, 5\}$, vertex 2 cannot move, while 4 and 5 do not satisfy $2E_i^Y > E_i^V$. Hence $T_{\text{noMovement}}$ becomes 2 and in the next iteration the algorithm terminates.*

*The resulting bipartite subgraph is determined by the bipartition $X = \{1, 3\}$ and $Y = \{2, 4, 5\}$, which gives $r^b = \frac{2}{3}$. Note that this ratio could be improved to $\frac{5}{6}$, if the initial partition was chosen so that it ends with $X = \{2, 3, 5\}$ and $Y = \{1, 4\}$.*

*2.2. Eigenvector-based methods*

In these methods, the initial partition is determined from the sign pattern of nonzero entries of $\nu^{\lambda_i^M}$, where $i$ corresponds to the smallest eigenvalue for $M = A$ and $M = Q$, and to the largest eigenvalue for $M = L$ and $M = \mathcal{L}$, while vertices with zero entries are randomly distributed to one of the parts. Afterwards, the movement routine of the local switching algorithm is applied in order to further increase the number of edges in the resulting bipartite subgraph.

*2.3. Identification of bipartite subgraphs using edge bipartivity*

Estrada and Rodríguez-Velázquez (Estrada and Rodríguez-Velázquez, 2005) quantified the degree of bipartivity of a complex network based on the ratio between the numbers of even and all closed walks:

$$\beta(G) = \frac{\sum_{j=1}^{V} \cosh(\lambda_j^A)}{\sum_{j=1}^{V} e^{\lambda_j^A}}. \tag{1}$$

The denominator above is the Estrada index (Estrada, 2000; de la Peña et al., 2007) which, through the expansion $e^x = \sum_{k \geq 0} \frac{x^k}{k!}$, results in the weighted series of all closed walks:

$$\sum_j e^{\lambda_j^A} = \sum_{k \geq 0} \frac{\text{Tr}(A^k)}{k!}.$$

Similarly, the numerator represents equally weighted series of all even walks:

$$\sum_j \cosh(\lambda_j^A) = \sum_{k \geq 0} \frac{\text{Tr}(A^{2k})}{(2k)!}.$$

It holds that $\frac{1}{2} < \beta(G) \leq 1$. If a graph is bipartite, it does not contain odd closed walks, so that $\beta(G) = 1$ in such case. The lower bound is reached in the limit $n \to \infty$ for the complete graph $G = K_n$.

Estrada and Rodríguez-Velázquez (Estrada and Gómez-Gardeñes, 2016) improved the lower bound to become 0 by introducing a new measure of bipartivity:

$$\beta(G)^{\text{new}} = \frac{\sum_{j=1}^{V} e^{-\lambda_j^A}}{\sum_{j=1}^{V} e^{\lambda_j^A}},$$

opening up the possibility to identify bipartite subgraphs with more edges. Nevertheless, both measures fail to distinguish among graphs cospectral with respect to $A$ (see Section 4 for discussion of an example). In this work, we take into account the improved measure $\beta(G)^{\text{new}}$.

The measure $\beta(G)^{\text{new}}$ is monotone in the sense that if addition of an edge $e$ increases the proportion of odd closed walks, then

$$\beta(G)^{\text{new}} \geq \beta(G+e)^{\text{new}}.$$

Hence contribution of a particular edge $e$ toward the bipartiteness of $G$ may be given by

$$\beta(e) = 1 - [\beta(G-e)^{\text{new}} - \beta(G)^{\text{new}}]. \tag{2}$$

A bipartite subgraph of $G$ may now be constructed by repeatedly removing an edge minimizing $\beta(e)$ until the graph becomes bipartite, as shown in Algorithm 2.

---

**Algorithm 2:** Identification of bipartite subgraphs using $\beta(e)$.

**Input**: $G = (V, E)$
**Result**: $r^b$

1 **while** do
2     **if** $G$ *is bipartite* **then**
3         compute $r^b$ ;
4         exit;
5     **else**
6         $e \leftarrow \arg\min_{e \in E(G)} \beta(e)$;
7         $G \leftarrow G - e$
8     **end**
9 **end**

---

Inspired by the edge bipartivity measure of Estrada and Rodríguez-Velázquez, we propose two new measures, denoted as $\Phi_A(e)$ and $\Phi_\mathcal{L}(e)$, based on the smallest eigenvector of $A$ and the largest eigenvector of $\mathcal{L}$, respectively. For an edge $e$ with the end vertices $i$ and $j$, $\Phi_A(e)$ is defined as:

$$\Phi_A(e) = \frac{\nu_i^{\lambda_n^A} \nu_j^{\lambda_n^A}}{\nu_i^{\lambda_1^A} \nu_j^{\lambda_1^A} + |\nu_i^{\lambda_n^A} \nu_j^{\lambda_n^A}|}, \tag{3}$$

where, as we recall, $\lambda_1^A > \lambda_2^A \geq \cdots \geq \lambda_n^A$. Since $G$ is assumed to be connected, by the Perron-Frobenius theorem (Perron, 1907; Pillai et al., 2005), $\lambda_1^A$ is a single eigenvalue with a positive eigenvector $\nu_i^{\lambda_1^A} > 0$ for $i \in V$. Addition of the absolute value of the term $\nu_i^{\lambda_n^A} \nu_j^{\lambda_n^A}$ in the denominator makes the value $\Phi_A(e)$ bounded as $-1 < \Phi_A(e) < 1$. In case of bipartite graph, each edge $e$ has $\Phi_A(e) = -\frac{1}{2}$, and edges that so-to-say stay in the way of bipartiteness are those with $\Phi_A(e) > 0$. Greedy approach then assumes that iterative removal of an edge with the maximum value of $\Phi_A(e)$ will make edge-deleted subgraph bipartite as quickly as possible.

The measure $\Phi_\mathcal{L}(e)$ may be defined similarly:

$$\Phi_\mathcal{L}(e) = \frac{\nu_i^{\lambda_1^\mathcal{L}} \nu_j^{\lambda_1^\mathcal{L}}}{\nu_i^{\lambda_n^\mathcal{L}} \nu_j^{\lambda_n^\mathcal{L}}}, \tag{4}$$

where $\lambda_1^{\mathcal{L}} > \lambda_2^{\mathcal{L}} \geq \cdots \geq \lambda_n^{\mathcal{L}}$. However as the eigenvector of the smallest eigenvalue of $\mathcal{L}$ is the all-one vector **1**, Eq. 4 reduces simply to:

$$\Phi_{\mathcal{L}}(e) = \nu_i^{\lambda_1^{\mathcal{L}}} \nu_j^{\lambda_1^{\mathcal{L}}}. \qquad (5)$$

Since $|\nu_i^{\lambda_1^{\mathcal{L}}} \nu_j^{\lambda_1^{\mathcal{L}}}| \leq \frac{\left(\nu_i^{\lambda_1^{\mathcal{L}}}\right)^2 + \left(\nu_j^{\lambda_1^{\mathcal{L}}}\right)^2}{2}$, the value $\Phi_{\mathcal{L}}(e)$ is bounded as well: assuming that $\max_i \nu_i^{\lambda_1^{\mathcal{L}}} = 1$ we get that $-1 \leq \Phi_{\mathcal{L}}(e) \leq 1$.

Similarly as in the case of $\beta(e)$, bipartite subgraphs of $G$ may be constructed by either repeatedly removing an edge that maximizes the value of $\Phi_A(e)$ or repeatedly removing an edge that maximizes the value of $\Phi_{\mathcal{L}}(e)$ until the graph becomes bipartite. Pseudo code for the method based on $\Phi_A(e)$ is shown in Algorithm 3. The method for $\Phi_{\mathcal{L}}(e)$ is obtained by replacing $A$ by $\mathcal{L}$.

---

**Algorithm 3:** Identification of bipartite subgraphs using $\Phi_A(e)$.

**Input**: $G = (V, E)$
**Result**: $r^b$

1 **while** do
2    **if** *G is bipartite* **then**
3      compute $r^b$ ;
4      exit;
5    **else**
6      $e \leftarrow \underset{e \in E(G)}{\arg\max} \dfrac{\nu_i^{\lambda_n^A} \nu_j^{\lambda_n^A}}{\nu_i^{\lambda_1^A} \nu_j^{\lambda_1^A} + |\nu_i^{\lambda_n^A} \nu_j^{\lambda_n^A}|}$ ;
7      $G \leftarrow (G - e)$
8    **end**
9 **end**

---

### 2.4. Parameter values for the graph models

For each of the four graph models: E-R, W-S, RG, and B-A, 1000 different graphs with 20 vertices are generated by uniformly sampling the respective parameter in the range specified in Table 2, and the ratio $r^b$ is calculated. To minimize the effect of random partitioning of vertices in the local switching algorithm, the maximum value of $r^b$ over 100 different random permutations (of vertices) are considered. For W-S graph model, the range is chosen such that the resulting graph will either be a regular graph (when $\psi$ is close to 0) or a graph with the small-world property because as $\psi$ approaches to 1, the graph tends to become an E-R type graph.

| Graph Models | Parameters | Sampled range [initial value final value] |
|---|---|---|
| E-R | $p$ | [0.2  1] |
| W-S | $(\psi, k)$ | ([0  0.3], 8) |
| RG | $(l, r)$ | (2, [0.5  1]) |
| B-A | $m$ | [1  10] |

**Table 2.** Parameter values for the graph models. For an E-A type graph, $p$ is the probability of attachment. For a W-S type graph, $\psi$ represents the probability of rewiring and $k$ is the mean degree. For a random geometric graph, $l$ represents the $\|.\|_l$ norm, while $r$ is the threshold for neighbor joining. For a B-A type graph, $m$ is the number of edges to attach in every step.

## 2.5. Software routines

All the algorithms described in this paper are implemented in Python. Generations of different graph models are carried out using the *NetworkX* (https://networkx.github.io/) package. For the purpose of reproducibility, codes are made available at https://github.com/DebCompBio/Eigenvector_based_bipartite_graphs.

## 3. Results

Figures 1, 2, 3, and 4 summarize the distribution of the $r^b$ values for the local switching and eigenvector-based methods on E-R, W-S, RG, and B-A graph models (see Appendix A for a brief overview) respectively. The distribution is obtained for 1000 different graphs corresponding to each graph model. Initially, we observe that for E-R and RG graph models and to some extent for W-S graph model as well, the $\beta(e)^{\text{new}}$ based method yields positive values of $\widehat{\text{pdf}}(r^b)$ much before other methods. This is undesirable as it means that there are graph instances for which $\beta(e)^{\text{new}}$ yields lower value of $r^b$; hence, smaller bipartite subgraphs, whereas other methods yield a higher value of $r^b$; hence, larger bipartite subgraphs. Following the same line of argument, left-skewed distributions are also undesirable as seen for the method based on the $L-$matrix for E-R (Fig. 1) and RG (Fig. 3) graph models.

In order to obtain a finer distinction between the methods, we proceed toward obtaining the empirical cumulative distribution (eCDF) function of the $r^b$ values for each graph model. The eCDF tells us about the probability that a randomly chosen $r^b$ value is less than a specified value. In this context, a higher value of the cumulative distribution function is undesirable as it indicates higher probability of bipartite subgraphs with lesser number of edges obtained using the respective method. Figure 5 and 6 represent the empirical cumulative distribution function (eCDF) of the $r^b$ values without and with $\beta(e)^{\text{new}}$, $\Phi_A(e)$, and $\Phi_\mathcal{L}(e)$ based methods respectively.

An even finer quantitative comparison between the methods in discussion has been made based on the fraction of graph instances for which a particular method yields higher and the same $r^b$ values to that of another method. The latter indicates the degree of similarity between two methods. Figures 7 and 8 show the heat-maps based on the fractions of graph instances representing superiority (*left* with *cool* color map) and similarity (*right* with *copper* color map) between methods. For example, consider the following two methods: $M_i$ and $M_j$, to obtain bipartite subgraphs. Additionally, consider $\text{Frac}(r^b_{M_{i>j(i=j)}})$ as the fraction of graph instances for which $r^b_{M_i}$ is higher than (same as) $r^b_{M_j}$. Now, the higher is the value of $\text{Frac}(r^b_{M_{i>j(i=j)}})$, as represented by the value of the cell $(i, j)$ in the heat-map, the better (more similar) is the method $M_i$ compare to (to) the method $M_j$.

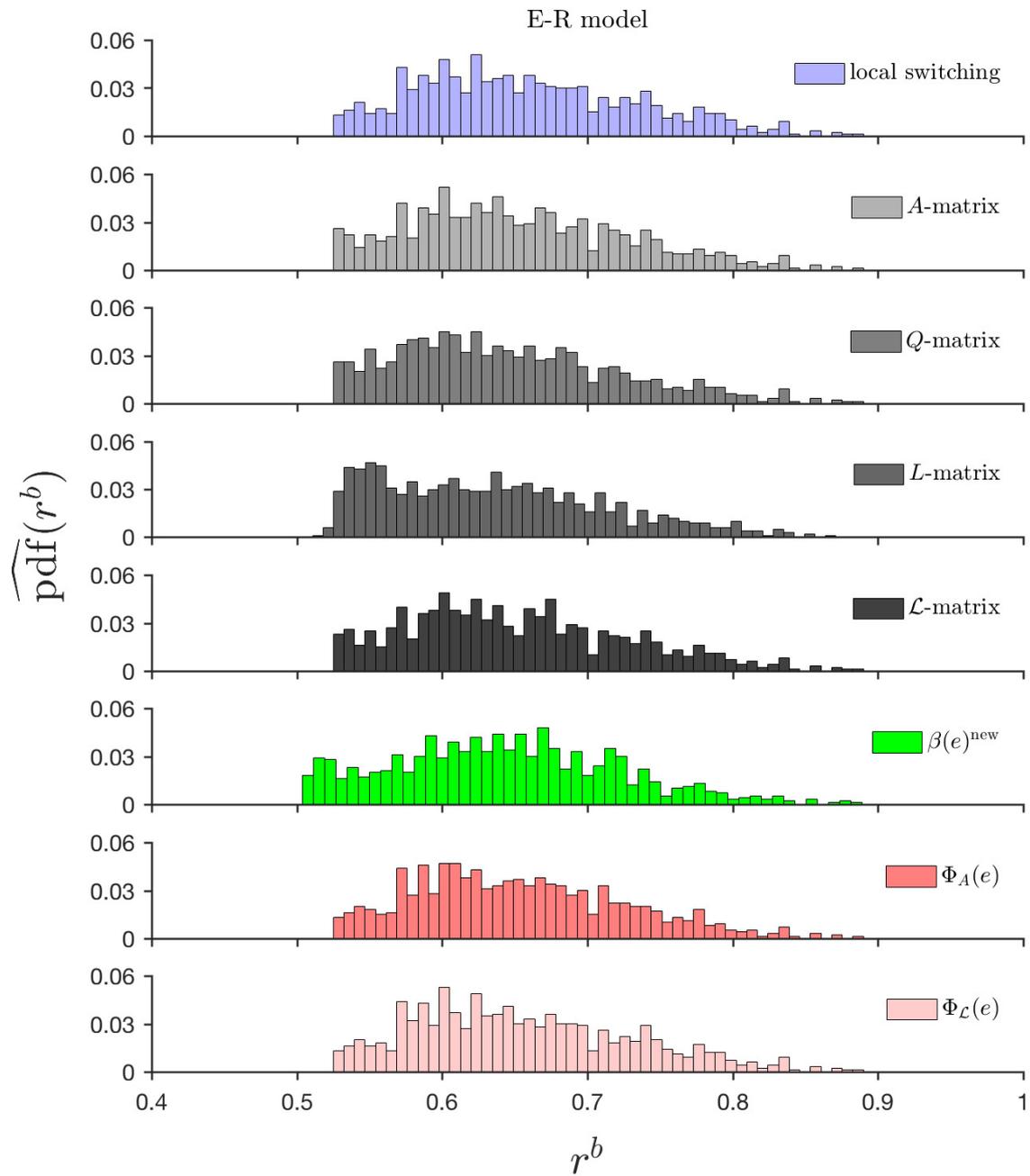

Fig. 1. Distribution of the $r^b$ values for the E-R graph model.

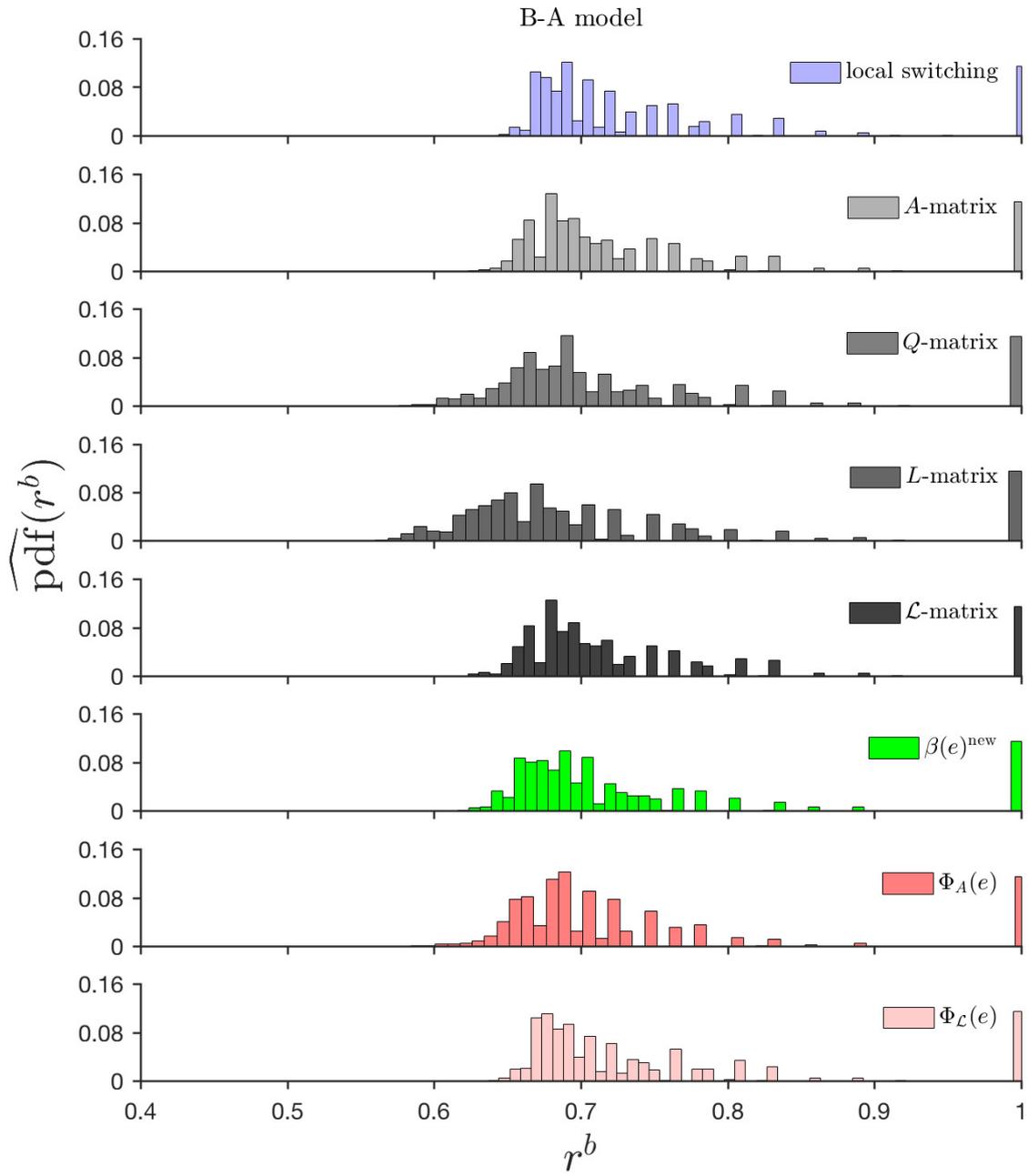

Fig. 2. Distribution of the $r^b$ values for the B-A graph model.

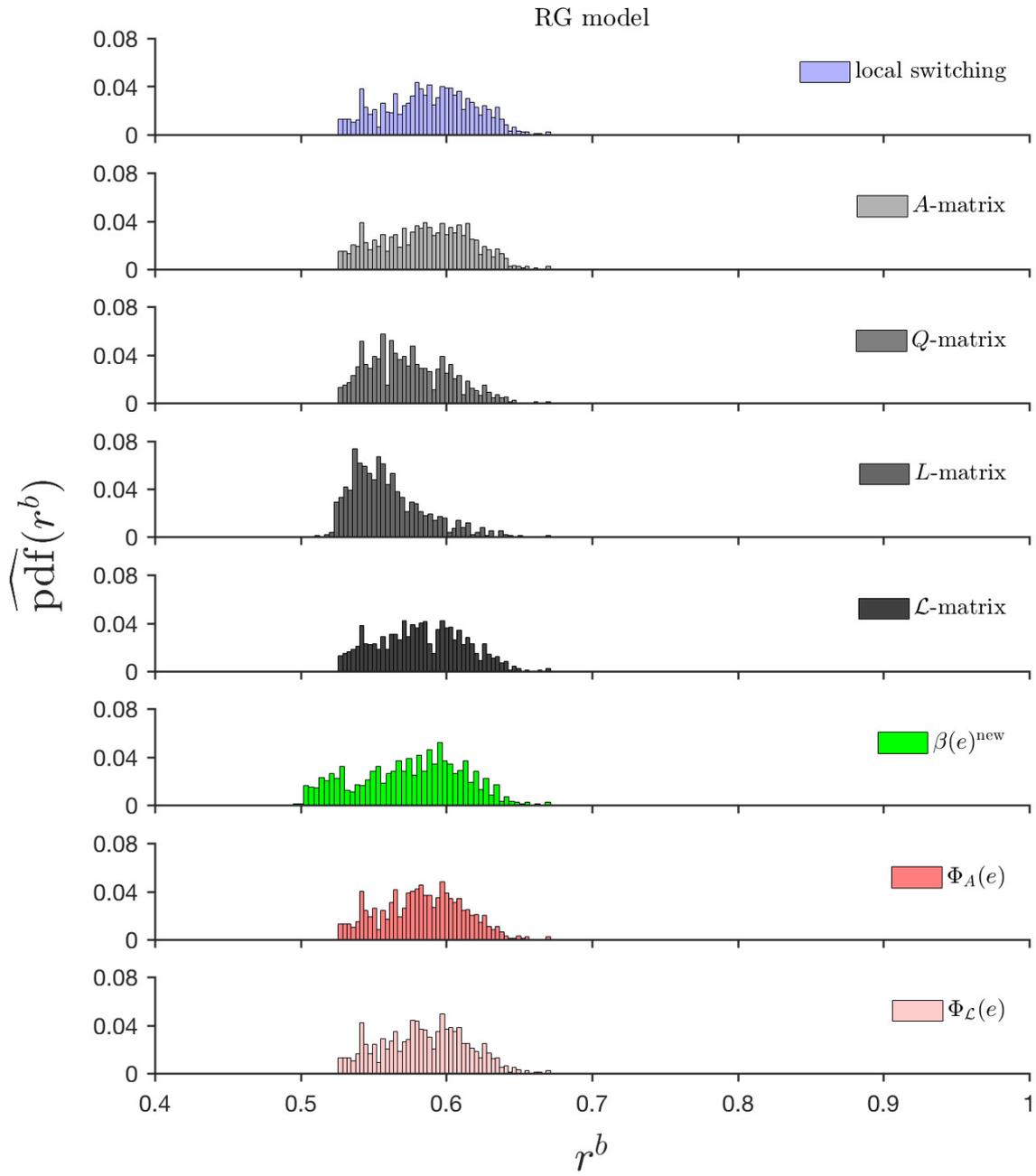

Fig. 3. Distribution of the $r^b$ values for the RG graph model.

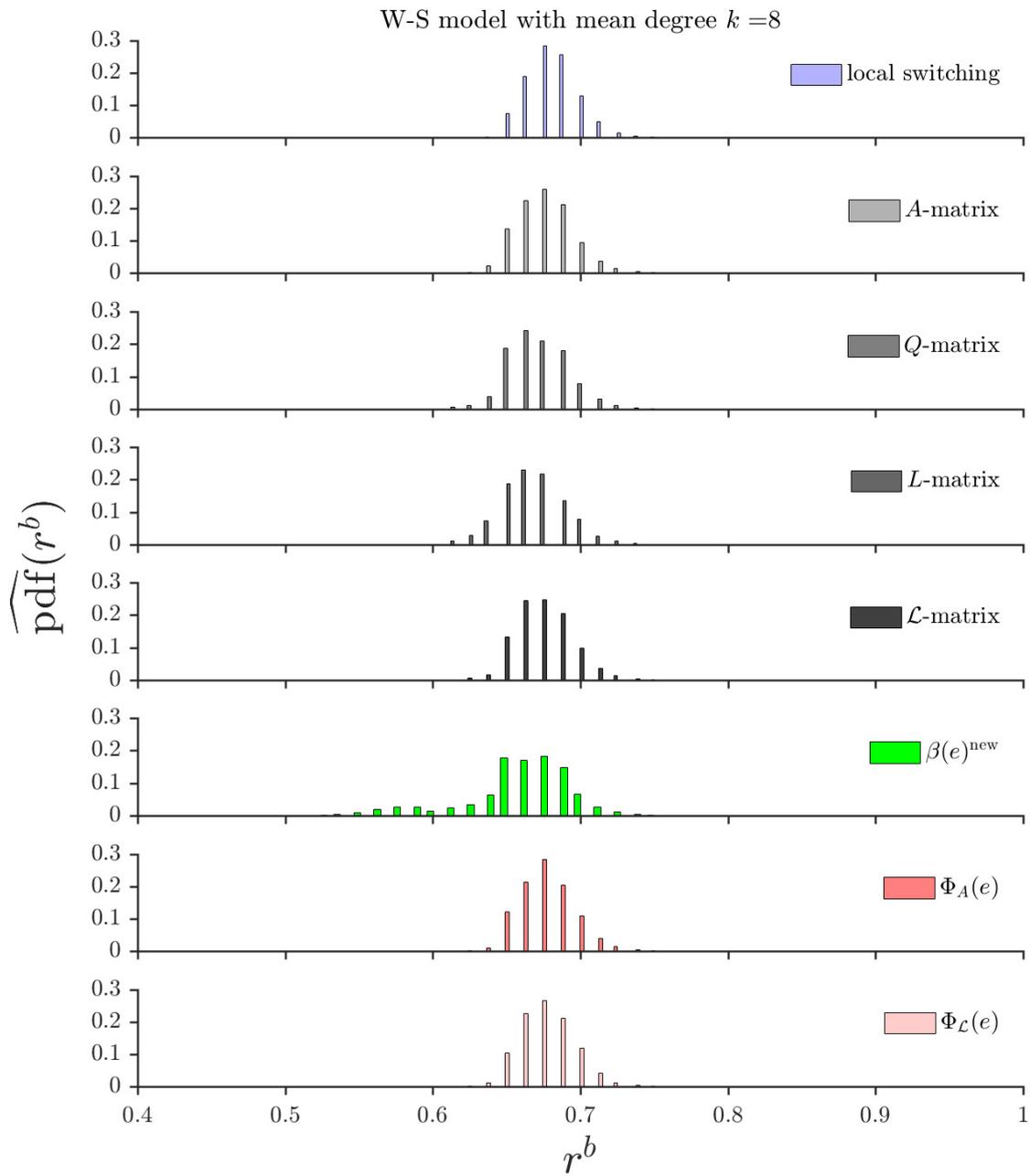

Fig. 4. Distribution of the $r^b$ values for the W-S graph model with mean degree 8.

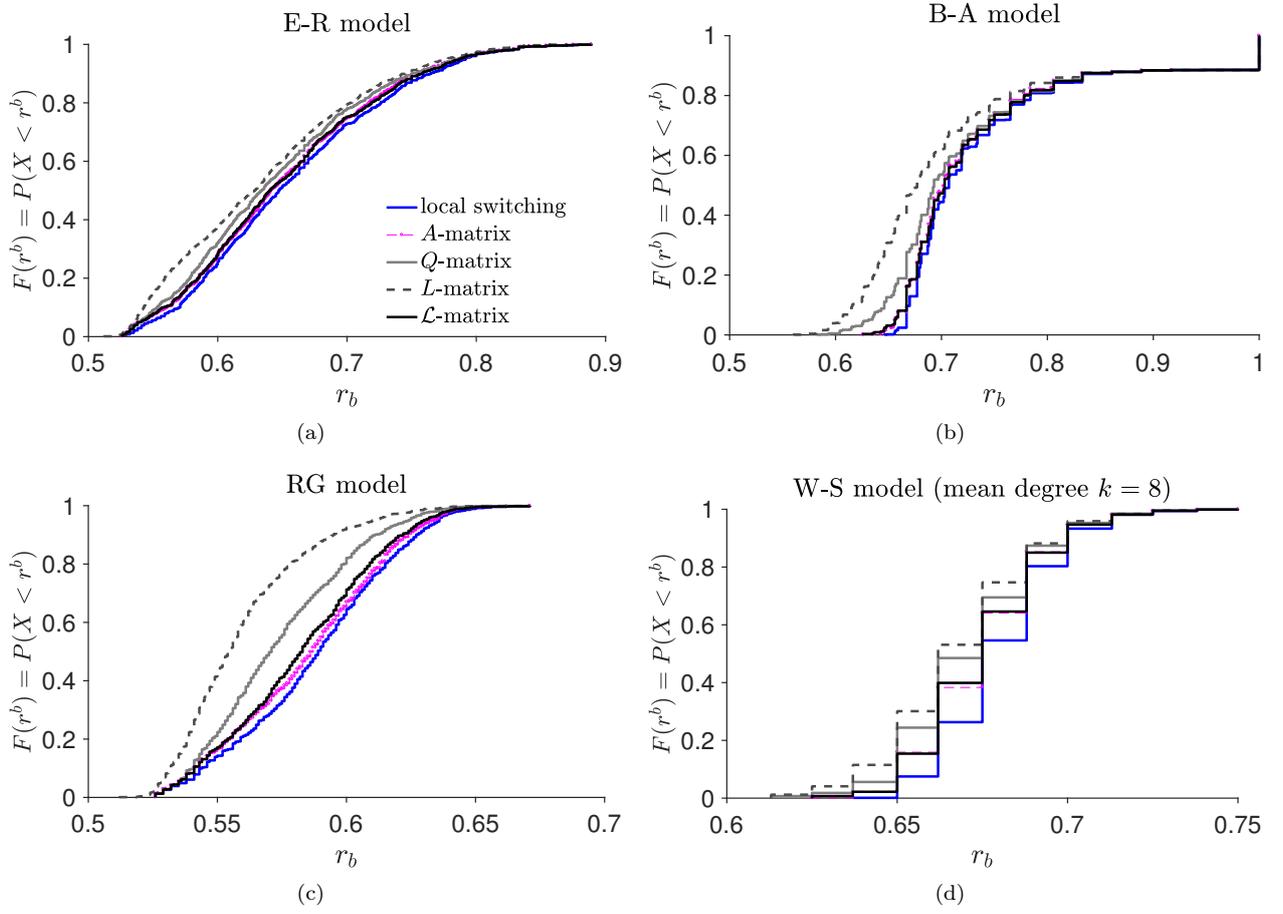

**Fig. 5. Empirical cumulative distribution of the $r^b$ values for local switching algorithm, $A, Q, L,$ and $\mathcal{L}$-matrix based methods.** Comparing local switching algorithm, $A, Q, L,$ and $\mathcal{L}$-matrix based methods for (a) E-R (b) B-A (c) RG (d) W-S graph models by constructing the empirical cumulative distributions corresponding to Figs 1, 2, 3, and 4 respectively.

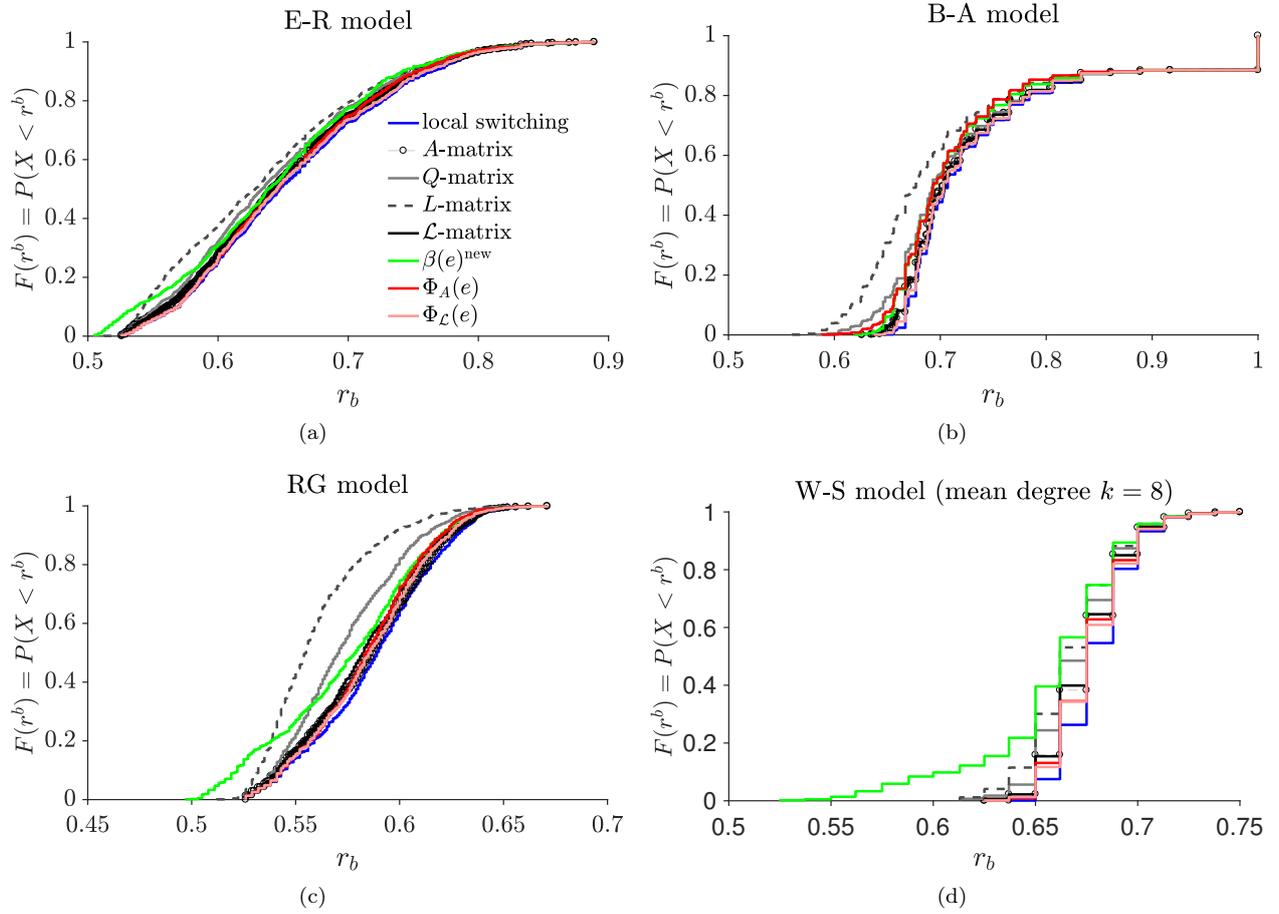

**Fig. 6. Empirical cumulative distribution of the $r^b$ values for** (a) E-R (b) B-A (c) RG (d) W-S graph models as in Fig. 5 but including $\beta(e)^{\text{new}}$, $\Phi_A(e)$, and $\Phi_\mathcal{L}(e)$ based methods.

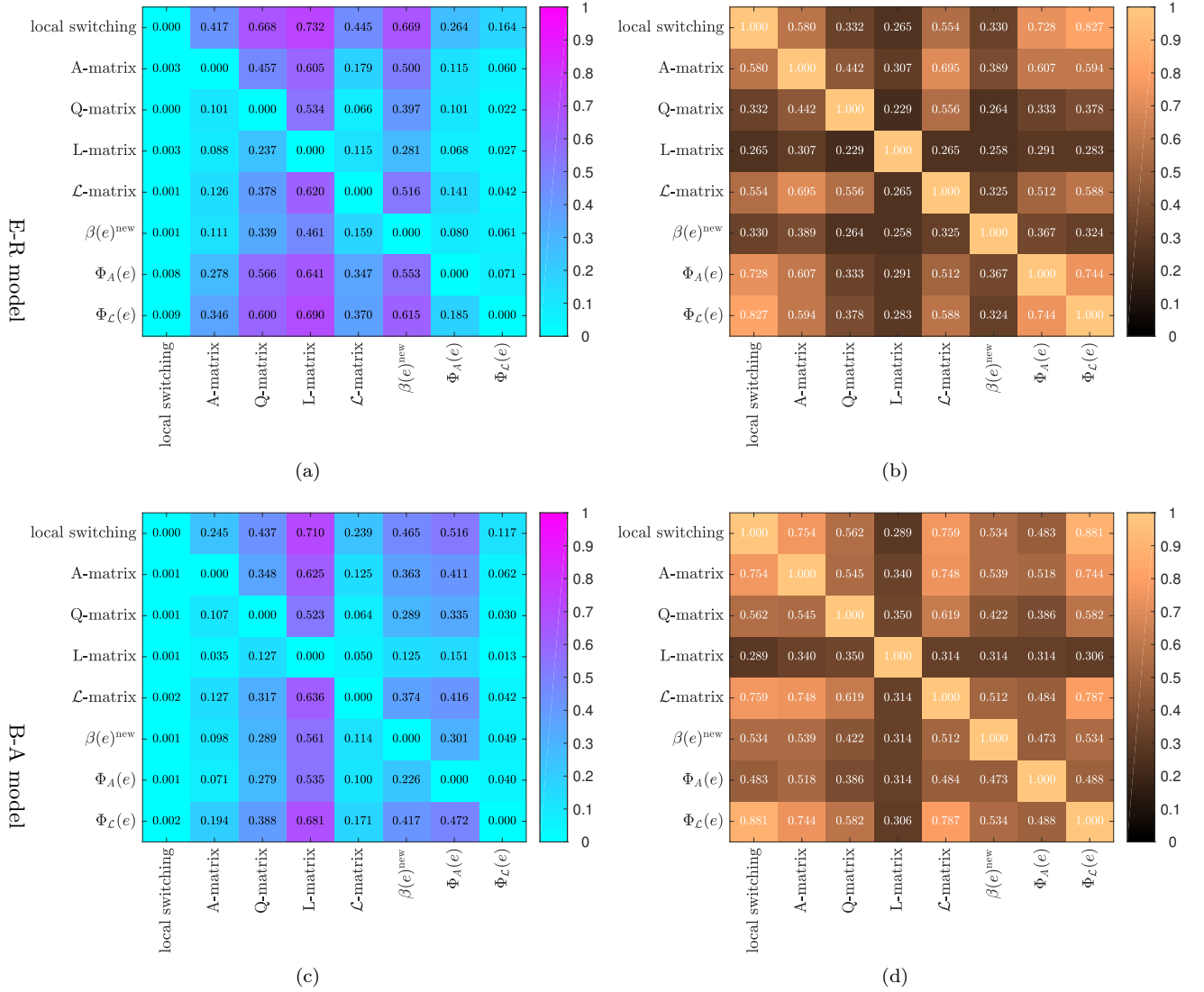

**Fig. 7. Heat-maps to compare different methods for E-R and B-A graph models.** (a)-(b) For E-R model, an entry $(i,j)$ in the heat-map matrix represents fraction of graph instances out of 1000 graphs when $r^b_{M_i} > r^b_{M_j}$ and $r^b_{M_i} = r^b_{M_j}$ respectively, where $r^b_{M_i}$ and $r^b_{M_j}$ corresponds to the value of $r^b$ for a method $M$ in the $i^{\text{th}}$ row and $j^{\text{th}}$ column respectively. (c)-(d) represent the same as (a)-(b) but for the B-A model.

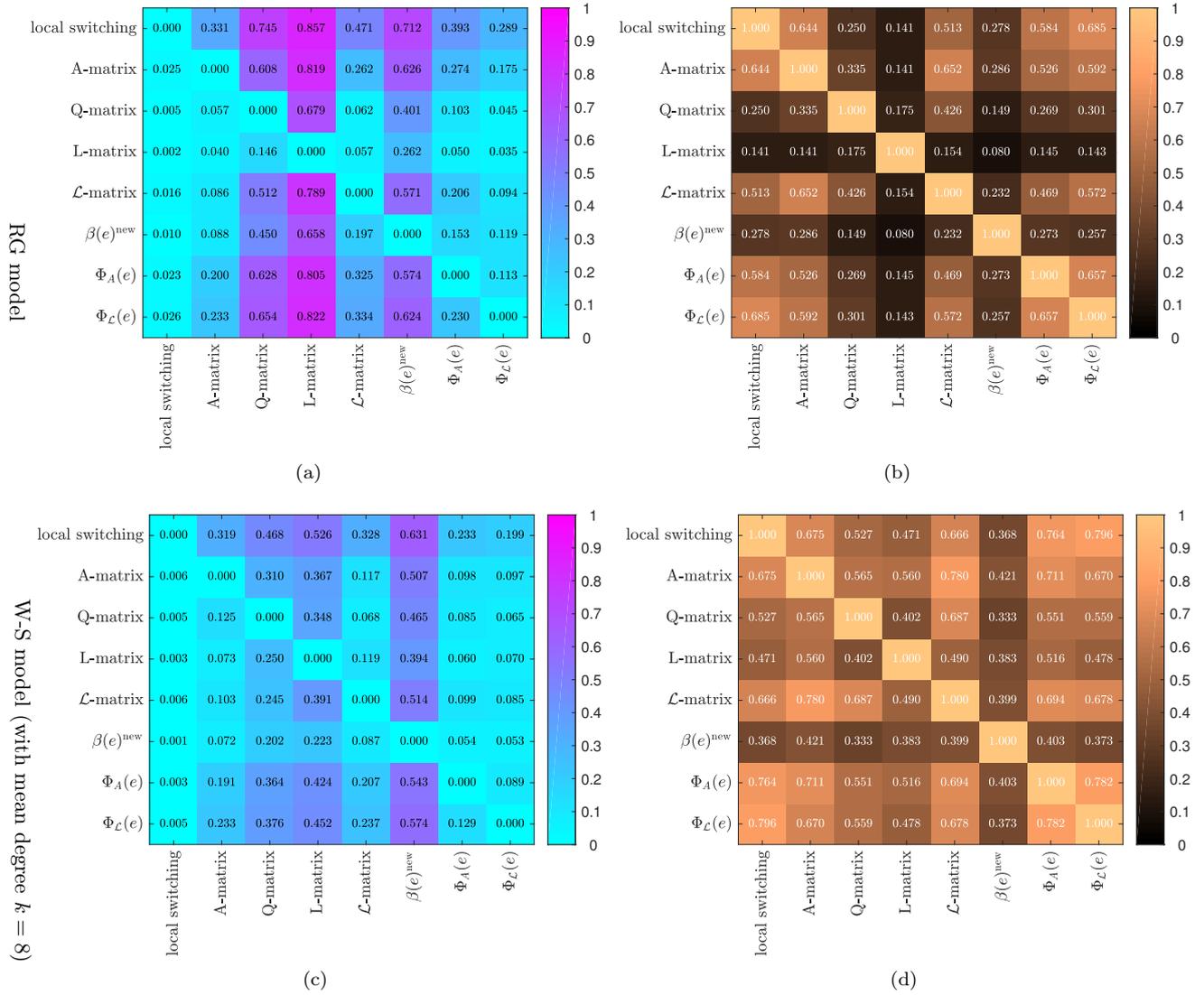

**Fig. 8. Heat-maps to compare different methods for RG and W-S graph models.** Same set of figures as in Fig. 7 but for RG (a)-(b), and for W-S (with mean degree $k = 8$) (c)-(d) graph models.

## 4. Discussion

In this work, we propose a set of eigenvector-based methods to identify bipartite subgraphs within a simple, undirected, connected, and non-trivial graph. At first, we compare the method based on sign-based partitioning to that of a local switching algorithm (see Algorithm 1) based on Erdös' bound. Subsequently, analogous to the notion of edge bipartivity index of Estrada and Gomez-Gardeñes denoted here by $\beta(e)^{\text{new}}$, we propose two new measures for edge bipartivity index, $\Phi_A(e)$ and $\Phi_\mathcal{L}(e)$ (see Eqs 3 and 4 respectively), based on the eigenvectors of $A$ and $\mathcal{L}$ matrices relying on their superior performances as evident from Figs (1 - 4) and Fig. 5. Experimental results over four different graph models: E-R, B-A, RG, and W-S, reveal that the local switching algorithm along with the $\Phi_A(e)$ and the $\Phi_\mathcal{L}(e)$ based methods outperform the rest. In fact, for graph models like RG and W-S, the local switching algorithm has been found to exhibit better performance than rest of the methods (see Figs 5(c) and 5(d)).

A finer distinction between the methods are illustrated in Figs 7 and 8 through heat-maps. The heat-maps quantitatively describe the degree of superiority (*left*) of a particular method over others as well as the degree of similarity (*right*) of that method with others. For example, in Fig. 7(a), for the E-R graph model, the cells in the row corresponding to the $L$-matrix mostly belongs to the lighter region of the spectrum in the color bar. It indicates that, out of 1000 graph instances, the $L$-matrix based method mostly yields lower $r^b$ values than the other methods for a significantly higher number (the number is represented as the fraction of the total number of graph instances) of graph instances and therefore towards its worst performance. It should be noted that Fig. 6(a) is in agreement with this observation as well. On the other hand, the respective rows corresponding to the local switching algorithm and the $\Phi_\mathcal{L}(e)$-based method contain comparatively higher number of cells that belong to the darker region of the spectrum indicating towards their superior performances. Next, while comparing the local switching algorithm and $\Phi_\mathcal{L}(e)$-based method, we observe that for around 0.9% of graph instances, the later yields higher $r^b$ values whereas for around 16.4% of graph instances, the former yields higher $r^b$ values; hence, performs slightly better. Next, according to the similarity matrix (Fig. 7(b)) for the E-R graph model, the local switching algorithm exhibits higher similarity with $A$-matrix, $\mathcal{L}$-matrix, $\Phi_A(e)$, and $\Phi_\mathcal{L}(e)$ based methods. In addition, one interesting observation is that the $\beta(e)^{\text{new}}$-based method neither exhibit a good performance nor exhibit good similarities to other methods. A concrete reasoning for this behavior is still lacking and can be posed as an open problem for the future.

In summary, we observe that the local switching algorithm outperforms rest of the methods for all the graph models. Moreover, as the runtime complexity of the switching algorithm is linearly dependent on the number of vertices, it is time-efficient for denser graphs as well unlike the other methods which are primarily edge-based. On the contrary, a potential drawback of the switching algorithm is that the update of the status of a vertex is based on its neighbors only, and because of this local updates the performance depends on the nature of the initial partition as illustrated in the Example 2.1. Therefore, one has to consider the maximum value over a number of different partition pairs. The optimum number of different partition pairs required for an arbitrary graph is still unknown; hence, can be posed as another open problem for the future as well. In this scenario, methods based on $\Phi_A(e)$ and $\Phi_\mathcal{L}(e)$ can be suggested as an alternative to the local switching algorithm that do not suffer from the problem based on the nature of initial partition, and in addition both the measures yield comparable results for all the graph models.

Although we propose two new edge bipartivity indices based on eigenvectors, we are unable to provide a measure for overall graph bipartivity unlike $\beta(G)^{\text{new}}$ of Estrada and Gomez-Gardeñes. In this context, it is worth noting that construction wise $\beta(G)^{\text{new}}$ relies on the $\{\lambda_i^A\}$ only, and therefore unable to distinguish between two co-spectral mates (with respect to $A$). For example, consider two co-spectral mates (non-isomorphic with same number of edges): $G_1$ and $G_2$ as depicted in the Figure 9. Additionally, consider that for $G_1$ it takes at least $e_1$ edges to be removed in order to make it bipartite. For $G_2$, it takes $e_1$ edges where $e_1 > e_2$. Therefore $G_2$ is *closer* than $G_1$ to become a bipartite graph which should be reflected as $\beta(G_2)^{\text{new}} > \beta(G_1)^{\text{new}}$. But, the values of $\beta^{\text{new}}(G)$ will remain same for both the graphs $G_1$ and $G_1$.

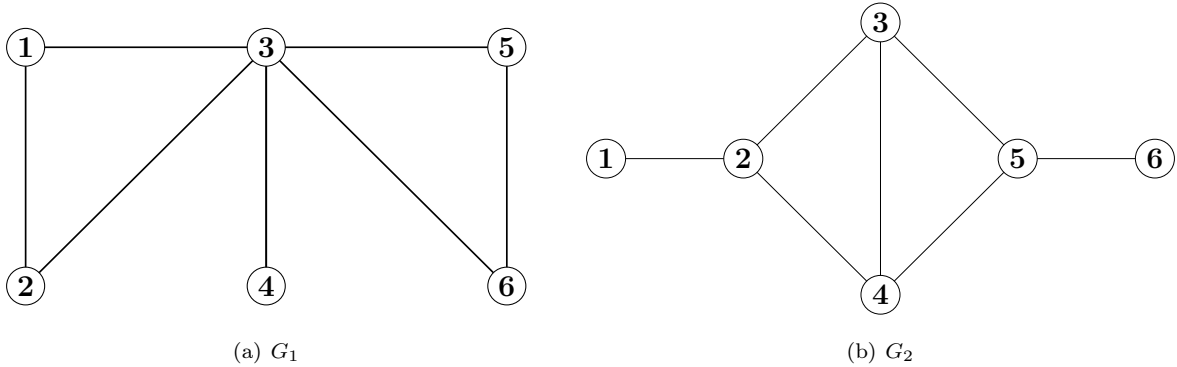

Fig. 9. Connected cospectral mates with respect to $A$ $G_1$ and $G_2$ have the same spectrum and number of edges. $G_1$ requires at least two edges (say $(1-2)$ and $(5-6)$, it can be some other pairs) to be removed to become bipartite while $G_2$ requires only one: $(3-4)$

Recall $\beta(G)^{new} = \frac{\sum_{j=1}^{V} e^{-\lambda_j^A}}{\sum_{j=1}^{V} e^{\lambda_j^A}}$. As $\{\lambda_i^A\}$ for $G_1$ and $G_2$ are same: $\begin{bmatrix} -1.9032 \\ -1.0000 \\ -1.0000 \\ 0.1939 \\ 1.0000 \\ 2.7093 \end{bmatrix}$, the values of $\beta(G)^{new}$ will be same for both the graphs. Now, while applying Algorithm 2 on both the graphs, $G_1$ requires at least two edges to be removed before it becomes bipartite, while $G_2$ requires only one edge:$(3-4)$ to be removed. The resulting bipartition will be $X = \{1, 3, 4, 6\}, Y = \{2, 5\}$. The same partition can be obtained from $G_1$ after removing edges $(1-3)$ and $(5-6)$. Therefore, $G_2$ is *more* bipartite or nearer to bipartite than $G_1$ but $\beta(G)^{new}$ is unable to distinguish that. The problem might be alleviated through a measure based on eigenvectors instead of eigenvalues but it depends on further properties whether they will also have the same eigenvector belonging to the minimal eigenvalue. This is possible without graphs being isomorphic but at the moment this is a problem worth to be addressed in future.

Overall, our study introduces eigenvector-based alternatives to the solution of the MAX-CUT problem which is, to the best of our knowledge, has not been addressed yet. Additionally, the study opens up new possibilities in addressing similar kind of problems based on the spectral graph theory.

## Appendix A. Graph models

*Erdös-Rényi graph model*

Erdös and Rényi in their 1959's paper (Erdös and Rényi, 1959) introduced the concept of a *random* graph model denoted by $G(V, E)$ with $|V|$ labeled vertices and $|E|$ edges. The graph is random in the sense that the configuration is chosen *randomly* from one of the $\binom{\binom{|V|}{2}}{|E|}$ configurations with equal probability. In another variant of the same graph model, denoted by $G(V, p)$, the nodes are connected randomly with probability $p$ (Erdös and Rényi, 1960). The connected edges are *independent* from each other. In this study, we use the $G(V, p)$ variant. For $G(V, p)$ to be almost surely connected, the condition is $p > \frac{(1+\epsilon)*\ln|V|}{|V|}$, for $0 < \epsilon < 1$.

*Random Geometric graph model*

The motivation for a Random Geometric graph lies in communication between two points in a $d$-dimensional space based on their *closeness*. In real life context, the points may be stations distributed across the whole country, connected nerve cells inside the brain cortex etc (Penrose, 2003). To put the description in a formal way, given a vertex set $V \in R^d$ and a norm $\|\cdot\|$ on $R^d$, where $R^d$ is a $d$-dimensional space on the set of real numbers $R$, the vertex pairs $\{(v_1, v_2) : v_1, v_2 \in V\}$ are connected with $\|v_1 - v_2\| \leq r$, where $r > 0$. Unlike the E-R graph model, the RG graph model does not follow the *property of independence* on edges, which means that a property defined on the pair of vertices of a RG graph is transitive in nature. This makes the RG graph model more *realistic* than that of a E-R graph model (Penrose, 2003).

*Watts-Strogatz graph model*

The Watts-Strogatz graph model addresses the *first* of the two shortcomings of the E-R graph models:

- Low clustering coefficient due to a constant, random and independent probability of connection between two nodes.
- Inability to account for power-law degree distribution which seems relevant for many biological, social, and technological networks (Watts and Strogatz, 1998).

A W-S graph model has a small *average path length* and high clustering coefficient. The construction of a W-S graph follows the two simple steps:

S-1: **Creation of ring lattice**: A ring lattice with $|V|$ nodes with mean degree $k$ is created where each node is connected to $k/2$ nearest neighbors.
S-2: **Rewiring of the target node**: For each edge in the graph, the target node is rewired with probability $\psi$

*Barabasi-Albert graph model*

The main concept of the Barabasi-Albert graph model is to *make a node richer (in links) if it is already rich*- a process called *preferential attachment* (Barabási and Albert, 1999). A node which is more connected is more likely to get new connections. This kind of preferential attachment alleviates the problem of not having a power law degree distribution which is more prevalent in real-world networks.

A shortcoming of the B-A model is that it fails to account for high-level clustering in real-world networks.

**Declarations**

*Funding*

This work is a part of a research collaboration between DP and DS initiated during a research visit of DP for three months at the University of Primorska, Koper, Slovenia. The reserach visit was funded by the Ministry of Education, Science, Culture and Sport of the Republic of Slovenia under bilateral mobility grant for the foreign nationals for 2010-11 (CMEPIUS).

*Author's contributions*

DP performed all the simulations. DS conceived the idea and and designed the study. DP wrote the main part of the manuscript. All authors discussed the results and implications and commented on the manuscript at all stages of the project. All authors read and approved the final manuscript.

*Declarations of interest*

None.

# References

Asratian AS, Denley TM, Häggkvist R. Bipartite graphs and their applications. volume 131. Cambridge University Press, 1998.

Baker E, Culpepper C, Philips C, Bubier J, Langston M, Chesler EJ. Identifying common components across biological network graphs using a bipartite data model. In: BMC proceedings. BioMed Central; volume 8; 2014. p. S4.

Barabási AL, Albert R. Emergence of scaling in random networks. Science 1999;286(5439):509–12.

Bylka S, Idzik A, Tuza Z. Maximum cuts: Improvements and local algorithmic analogues of the edwards-erdos inequality. Discrete Mathematics 1999;194(1-3):39–58.

Chao C, Thomaz AL. Timing in multimodal turn-taking interactions: Control and analysis using timed petri nets. Journal of Human-Robot Interaction 2012;1(1).

Chung FR. Spectral graph theory, CBMS Regional Conference Series in Mathematics, No. 92. American Mathematical Society, 1996.

Dall J, Christensen M. Random geometric graphs. Physical Review E 2002;66(1):016121.

Edwards C. Some extremal properties of bipartite subgraphs. Canadian Journal of Mathematics 1973;25(3):475–83.

Erdös P. On some extremal problems in graph theory. Israel Journal of Mathematics 1965;3(2):113–6.

Erdös P, Rényi A. On random graphs i. Publicationes Mathematicae Debrecen 1959;6:290–7.

Erdös P, Rényi A. On the evolution of random graphs. Publications of the Mathematical Institute of the Hungarian Academy of Science 1960;5:17–61.

Estrada E. Characterization of 3d molecular structure. Chemical Physics Letters 2000;319(5-6):713–8.

Estrada E, Gómez-Gardeñes J. Network bipartivity and the transportation efficiency of european passenger airlines. Physica D: Nonlinear Phenomena 2016;323:57–63.

Estrada E, Rodríguez-Velázquez JA. Spectral measures of bipartivity in complex networks. Physical Review E 2005;72(4):046105.

Gilbert EN. Random plane networks. Journal of the Society for Industrial and Applied Mathematics 1961;9(4):533–43.

Goemans MX, Williamson DP. . 879-approximation algorithms for max cut and max 2sat. In: Proceedings of the twenty-sixth annual ACM symposium on Theory of computing. ACM; 1994. p. 422–31.

Grone R, Merris R, Sunder VS. The laplacian spectrum of a graph. SIAM Journal on Matrix Analysis and Applications 1990;11(2):218–38.

Guillaume JL, Latapy M. Bipartite structure of all complex networks. Information processing letters 2004;90(5):215–21.

Guillaume JL, Latapy M. Bipartite graphs as models of complex networks. Physica A: Statistical Mechanics and its Applications 2006;371(2):795–813.

Håstad J. Some optimal inapproximability results. Journal of the ACM (JACM) 2001;48(4):798–859.

Holme P, Liljeros F, Edling CR, Kim BJ. Network bipartivity. Physical Review E 2003;68(5):056107.

Michael RG, David SJ. Computers and intractability: a guide to the theory of np-completeness. WH Free Co, San Fr 1979;:90–1.

Mitzenmacher M, Upfal E. Probability and computing: Randomized algorithms and probabilistic analysis. Cambridge University Press, 2005.

Moon J, Friedberg I, Eulenstein O. Highly bi-connected subgraphs for computational protein function annotation. In: International Computing and Combinatorics Conference. Springer; 2016. p. 573–84.

Moon TK. Error correction coding: mathematical methods and algorithms. Wiley, 2005.

Motwani R. Randomized algorithms. Cambridge University Press, 1995.

Newman ME. 2 random graphs as models of networks. Handbook of graphs and networks 2003;:35.

Newman ME, Watts DJ, Strogatz SH. Random graph models of social networks. Proceedings of the National Academy of Sciences of the United States of America 2002;99(Suppl 1):2566–72.

Pavlopoulos GA, Kontou PI, Pavlopoulou A, Bouyioukos C, Markou E, Bagos PG. Bipartite graphs in systems biology and medicine: a survey of methods and applications. GigaScience 2018;.

Pavlopoulos GA, Secrier M, Moschopoulos CN, Soldatos TG, Kossida S, Aerts J, Schneider R, Bagos PG. Using graph theory to analyze biological networks. BioData mining 2011;4(1):1.

de la Peña JA, Gutman I, Rada J. Estimating the estrada index. Linear Algebra and its Applications 2007;427(1):70–6.

Penrose M. Random geometric graphs. volume 5. Oxford University Press, 2003.

Perron O. Grundlagen für eine theorie des jacobischen kettenbruchalgorithmus. Mathematische Annalen 1907;64(1):1–76.

Pillai SU, Suel T, Cha S. The perron-frobenius theorem: some of its applications. IEEE Signal Processing Magazine 2005;22(2):62–75.

Platig J, Castaldi PJ, DeMeo D, Quackenbush J. Bipartite community structure of eqtls. PLoS Computational Biology 2016;12(9):e1005033.

Pržulj N, Corneil DG, Jurisica I. Modeling interactome: scale-free or geometric? Bioinformatics 2004;20(18):3508–15.

Roth R. On the eigenvectors belonging to the minimum eigenvalue of an essentially nonnegative symmetric matrix with bipartite graph. Linear Algebra and its Applications 1989;118:1–10.

Roy S, Edan Y. Investigating joint-action in short-cycle repetitive handover tasks: The role of giver versus receiver and its implications for human-robot collaborative system design. International Journal of Social Robotics 2018;:1–16.

Sachs H. Beziehungen zwischen den in einem graphen enthaltenen kreisen und seinem charakteristischen polynom. Publicationes Mathematicae Debrecen 1964;11:119–34.

Schweitzer F, Fagiolo G, Sornette D, Vega-Redondo F, Vespignani A, White DR. Economic networks: The new challenges. Science 2009;325(5939):422.

Watts DJ, Strogatz SH. Collective dynamics of 'small-world' networks. Nature 1998;393(6684):440–2.